# Constructing Dynamical Systems to Model Higher Order Ising Spin Interactions and their Application in Solving Combinatorial Optimization Problems

Mohammad Khairul Bashar, Nikhil Shukla*

Department of Electrical and Computer Engineering, University of Virginia,

Charlottesville, VA- 22904 USA

*e-mail: ns6pf@virginia.edu

**ABSTRACT**

The Ising model provides a natural mapping for many computationally hard combinatorial optimization problems (COPs). Consequently, dynamical system-inspired computing models and hardware platforms that minimize the Ising Hamiltonian, have recently been proposed as a potential candidate for solving COPs, with the promise of significant performance benefit. However, the Ising model, and consequently, the corresponding dynamical system-based computational models primarily consider quadratic interactions among the nodes. Computational models considering higher order interactions among Ising spins remain largely unexplored. Therefore, in this work, we propose dynamical-system-based computational models to consider higher order (>2) interactions among the Ising spins, which subsequently, enables us to propose computational models to directly solve many COPs that entail such higher order interactions (COPs on hypergraphs). Specifically, we demonstrate our approach by developing dynamical systems to compute the solution for the Boolean NAE-K-SAT (K≥4) problem as well as solve the Max-K-Cut of a hypergraph. Our work advances the potential of the physics-inspired 'toolbox' for solving COPs.

The minimization of the Ising Hamiltonian using dynamical systems such as coupled electronic [1]-[5] and photonic oscillators [6],[8] has received substantial attention in recent years [9],[10]. A significant driving force behind the effort to realize a so-called 'Ising machine' is that the solution to the Ising model can be mapped to many computationally intractable problems in combinatorial optimization (e.g., MaxCut, Traveling Salesman Problem (TSP) among others) [11]-[17]. Consequently, this creates the possibility of realizing Ising machine-inspired custom accelerators that can offer the possibility of significant performance benefits. However, dynamical system formulations that have been used to 'solve' the Ising model, typically consider only pair-wise coupling. From an application standpoint, while these characteristics capture quadratic interactions, the dynamical systems and their supporting computational models cannot be applied *directly* to solve problems that require higher order interaction among the spins [18],[19]. Therefore, the objective of this work is two-fold: (1) define dynamical systems that model higher order (>2) interactions among the Ising spins; and (2) map the resulting dynamics to relevant computational problems. We consider two examples: computing the solutions for the NAE-K-SAT (Not-All-Equal SAT) problem and computing the Max-K-Cut of a hypergraph. We emphasize here that our focus is on defining the system dynamics that capture the higher order interactions, and presently not on the physical implementation of the higher order interactions.

The general form to represent higher order interactions among the Ising spins can be expressed as,

$$H = -\sum_{i,j} J_{ij}^{(2)} s_i s_j - \sum_{i,j,k} J_{ijk}^{(3)} s_i s_j s_k - \sum_{i,j,k,l} J_{ijkl}^{(4)} s_i s_j s_k s_l \dots \quad (1)$$

Where $J_{ij}^{(2)}$ represents the pairwise interaction coefficient between two Ising spins. The first term on the right-hand side $(-\sum_{i,j} J_{ij}^{(2)} s_i s_j)$ is usually considered when describing quadratic/pairwise interactions among Ising spins $s = \{-1,1\}^n$. However, considering the higher order interactions among the spins can help describe the objective functions of several combinatorial optimization problems (COPs) as illustrated here with the example of the NAE-K-SAT problem (without the need for problem decomposition). The NAE-K-SAT problem is a constrained version of the SAT problem where the objective is to find an assignment for the variables of the given Boolean expression (in the conjunctive normal form) such that: (a) at least one variable in every clause is TRUE (i.e., the clause is satisfied; standard SAT constraint); (b) at least one variable in every clause is FALSE [20]. Using an approach similar to SAT, the NAE-K-SAT problem can be expressed as computing an assignment for the variables such that $Y (= C_1.S_1 \wedge C_2.S_2 \wedge \dots \wedge C_M.S_M) = 1$. Here, $C_i \equiv (x_1 \vee x_2 \vee \bar{x}_3 \dots \bar{x}_N)$, and $S_i \equiv (\bar{x}_1 \vee \bar{x}_2 \vee x_3 \dots x_N)$ (i.e., $S_i$ and $C_i$ have the same variables but in opposite forms). Traditionally, when considering only pairwise interactions among the Ising spins, mapping such problems entails significant pre-processing including the use of auxiliary variables that can significantly increase the size of the problem that must be eventually solved [19],[21]-[24] using the dynamical system.

**NAE-4-SAT:** To illustrate how we can map the NAE-K-SAT problem to higher order interactions among the Ising spins, we first consider the example of the NAE-4-SAT problem where each clause of the NAE-4-SAT problem consists of 4 literals, expressed in the general form as $(x_i \vee x_j \vee x_k \vee x_l).(\bar{x}_i \vee \bar{x}_j \vee \bar{x}_k \vee \bar{x}_l) \equiv (x_i \oplus x_j) \vee (x_i \oplus x_k) \vee$ $(x_i \oplus x_l) \vee (x_j \oplus x_k) \vee (x_j \oplus x_l) \vee (x_k \oplus x_l)$, where $x \in \{0,1\}^n$ ($x$ is a set of Boolean variables). K=4 is specifically chosen since it is the lowest K where higher order interactions among the Ising spins are required to formulate the objective function for the problem (shown in Table 1). To formulate the problem in terms of Ising spins, we utilize the following property among the Boolean variables and the spins $(x_i \oplus x_j) \equiv \frac{1-s_i s_j}{2}$. Here, the logic level 0 (1) corresponds to an evaluation of -1(1) of the expression on the right-hand side, respectively. Furthermore, the complement of the logical OR among the XOR terms $((x_i \oplus x_j) \vee (x_i \oplus x_k) \vee \ldots \vee (x_k \oplus x_l))$ can be expressed as, $\left(1 - \left(\frac{1-s_i s_j}{2}\right)\right) \cdot \left(1 - \left(\frac{1-s_i s_k}{2}\right)\right) \ldots \left(1 - \left(\frac{1-s_k s_l}{2}\right)\right)$. Simplifying the above expression yields $\left(\frac{1+s_i s_j}{2}\right)\left(\frac{1+s_i s_k}{2}\right)\left(\frac{1+s_i s_l}{2}\right) \ldots \left(\frac{1+s_k s_l}{2}\right) \equiv \frac{1}{8}(1 + s_i s_j + s_i s_k + s_i s_l + s_j s_k + s_j s_l + s_k s_l + s_i s_j s_k s_l)$.

It can be observed that besides the second order interaction terms, the resulting expression also contains a 4$^{th}$ order interaction term among the spins. Consequently, the objective function for the NAE-4-SAT problem, over M clauses, can be formulated as the minimization of

$$H_{NAE-4-SAT} = -\sum_{m=1}^{M} \left( \sum_{\substack{i,j \\ i<j}}^{N} (-c_{mi} c_{mj} s_i s_j) + \sum_{\substack{i,j,k,l \\ i<j<k<l}}^{N} (-c_{mi} c_{mj} c_{mk} c_{ml} s_i s_j s_k s_l) \right) \quad (2)$$

Here, $c_{mi} = 1(-1)$, if the $i^{th}$ variable appears in the $m^{th}$ clause in the normal (negated) form; $c_{mi} = 0$ if the $i^{th}$ variable is absent from the $m^{th}$ clause. Using the same approach, we derive such expressions for a few other values of K in the NAE-K-SAT problem in Table 1. Details of the derivation of the objective function for NAE-5-SAT are shown in Appendix I.

| K | Expression for a single clause & objective function for the NAE-K-SAT |
|---|---|
| 2 | **Expression for a single clause:** $$(x_i \vee x_j).(\bar{x}_i \vee \bar{x}_j) \equiv s_i s_j$$ **Objective function:** $$H = -\sum_{m=1}^{M} \sum_{i,j,i<j}^{N} (-c_{mi}c_{mj}s_i s_j) \equiv -\sum_{m=1}^{M} \sum_{i,j,i<j}^{N} J_{ij} s_i s_j$$ Where $J_{ij} = -c_{mi}c_{mj}$. *It can be observed that when the variables appear only in the normal form i.e., $c_{mi} \geq 0$, the expression represents the solution to the archetypal MaxCut problem.* |
| 3 | **Expression for a single clause:** $$(x_i \vee x_j \vee x_k).(\bar{x}_i \vee \bar{x}_j \vee \bar{x}_k) \equiv s_i s_j + s_i s_k + s_j s_k$$ **Objective function:** $$H = -\sum_{m=1}^{M} \sum_{i,j,i<j}^{N} (-c_{mi}c_{mj}s_i s_j)$$ |
| 4 | **Expression for a single clause:** $$(x_i \vee x_j \vee x_k \vee x_l).(\bar{x}_i \vee \bar{x}_j \vee \bar{x}_k \vee \bar{x}_l)$$ $$\equiv s_i s_j + s_i s_k + s_i s_l + s_j s_k + s_j s_l + s_k s_l + s_i s_j s_k s_l$$ **Objective function:** $$H = -\sum_{m=1}^{M} \left( \sum_{\substack{i,j \\ i<j}}^{N} (-c_{mi}c_{mj}s_i s_j) + \sum_{\substack{i,j,k,l \\ i<j<k<l}}^{N} (-c_{mi}c_{mj}c_{mk}c_{ml}s_i s_j s_k s_l) \right)$$ |
| 5 | **Expression for a single clause:** $$(x_i \vee x_j \vee x_k \vee x_l \vee x_m).(\bar{x}_i \vee \bar{x}_j \vee \bar{x}_k \vee \bar{x}_l \vee \bar{x}_m)$$ $$\equiv s_i s_j + s_i s_k + s_i s_l + s_i s_m + s_j s_k + s_j s_l + s_j s_m + s_k s_l + s_k s_m + s_l s_m + s_i s_j s_k s_l$$ $$+ s_i s_j s_k s_m + s_i s_j s_l s_m + s_i s_k s_l s_m + s_j s_k s_l s_m$$ **Objective function:** |

$$H = -\sum_{m=1}^{M}\left(\sum_{\substack{i,j \\ i<j}}^{N}(-c_{mi}c_{mj}s_is_j) + \sum_{\substack{i,j,k,l \\ i<j<k<l}}^{N}(-c_{mi}c_{mj}c_{mk}c_{ml}s_is_js_ks_l)\right)$$

*We note that constants and scalars have not been shown here in the expression for the single clause as well as for the objective function.*

*Table 1.* *Objective functions for the NAE-K-SAT problem expressed using Ising spins*

**Constructing a dynamical system for the NAE-K-SAT problem:** We now aim to formulate the dynamical system and the corresponding energy function for the NAE-K-SAT problem; the dynamical system is designed such that the ground state of the 'energy' function (more precisely, the Lyapunov function) must correspond to a global optimum of the objective function. To construct this system, we draw inspiration from the dynamics of coupled oscillators under second harmonic injection which effectively forces the oscillator states to assume a binary phase value of 0 or π (details of the second harmonic injection can be found in work by Wang *et al* [25]). Without loss of generality, we assume that one spin state (say, $s = +1$) is represented by phase 0 while the other spin state ($s = -1$) is represented by the phase angle π. Subsequently, the second order interaction terms among the Ising spins $s_is_j$ can be represented by $\cos(\phi_i - \phi_j)$. When the spins are in opposite states i.e., $s_i = 1(-1); s_j = -1(1)$, $s_js_j \equiv \cos(\phi_i - \phi_j) = -1$, whereas when the spins are in the same states i.e., $s_i = 1(-1); s_j = 1(-1)$, $s_js_j \equiv \cos(\phi_i - \phi_j) = 1$. Similarly, the higher order interactions of even order can be modeled as shown in Table 2.

| Order | Ising interaction | Equivalent formulation for constructing dynamical system |
|---|---|---|
| 2 | $s_i s_j$ | $\cos(\phi_i - \phi_j)$ |
| 4 | $s_i s_j s_k s_l$ | $\cos(\phi_i - \phi_j + \phi_k - \phi_l)$ |
| 6 | $s_i s_j s_k s_l s_m s_n$ | $\cos(\phi_i - \phi_j + \phi_k - \phi_l + \phi_m - \phi_n)$ |

*Table 2.* Equivalent energy function for modeling higher order interactions among Ising spins. The second harmonic signal included as a part of the dynamics (not shown here) helps force $\phi$ to $\{0, \pi\}$.

The equivalence between the higher order terms and the corresponding energy term is shown in Table 3.

| <span style="color:red">**Second Order Interactions ($s_i . s_j$)**</span> | | | |
|---|---|---|---|
| $s_i \, s_j$ | $s_i . s_j$ | $\phi_i \, \phi_j$ | $\cos(\phi_i - \phi_j)$ |
| -1 -1 | +1 | π π | +1 |
| -1 +1 | -1 | π 0 | -1 |
| +1 -1 | -1 | 0 π | -1 |
| +1 +1 | +1 | 0 0 | +1 |
| <span style="color:red">**Fourth Order Interactions ($s_i . s_j . s_k . s_l$)**</span> | | | |
| $s_i \, s_j \, s_k \, s_l$ | $s_i . s_j . s_k . s_l$ | $\phi_i \, \phi_j \, \phi_k \, \phi_l$ | $\cos(\phi_i - \phi_j + \phi_k - \phi_l)$ |
| -1 -1 -1 -1 | +1 | π π π π | +1 |
| -1 -1 -1 +1 | -1 | π π π 0 | -1 |
| -1 -1 +1 -1 | -1 | π π 0 π | -1 |
| -1 -1 +1 +1 | +1 | π π 0 0 | +1 |
| -1 +1 -1 -1 | -1 | π 0 π π | -1 |
| -1 +1 -1 +1 | +1 | π 0 π 0 | +1 |
| -1 +1 +1 -1 | +1 | π 0 0 π | +1 |
| -1 +1 +1 +1 | -1 | π 0 0 0 | -1 |
| +1 -1 -1 -1 | -1 | 0 π π π | -1 |
| +1 -1 -1 +1 | +1 | 0 π π 0 | +1 |
| +1 -1 +1 -1 | +1 | 0 π 0 π | +1 |
| +1 -1 +1 +1 | -1 | 0 π 0 0 | -1 |
| +1 +1 -1 -1 | +1 | 0 0 π π | +1 |
| +1 +1 -1 +1 | -1 | 0 0 π 0 | -1 |
| +1 +1 +1 -1 | -1 | 0 0 0 π | -1 |
| +1 +1 +1 +1 | +1 | 0 0 0 0 | +1 |

***Table 3.*** *Equivalence between the higher order Ising spin interaction terms and the constructed energy functions.*

Using the above relationships developed in Table 1, the energy functions for the NAE-K-SAT problem can be formulated as shown in Table 4. The corresponding dynamics for which the above function is a Lyapunov function is also shown in Table 4.

| K | Objective function, equivalent energy function, and dynamics. |
|---|---|
| 2 & 3 | **Objective function:** $$H = -\sum_{m=1}^{M} \sum_{i,j,i<j}^{N} \left(-c_{mi}c_{mj}s_i s_j\right)$$ **Energy function:** $$E = C \sum_{m=1}^{M} \left[\sum_{i,j,i<j}^{N} c_{mi}c_{mj} \cos(\phi_i - \phi_j) + 1\right] - \frac{C_s}{2} \sum_{i=1}^{N} \cos(2\phi_i)$$ **Dynamics:** $$\frac{d\phi_i}{dt} = C \left[\sum_{m=1}^{M} \sum_{j=1}^{N} c_{mi}c_{mj} \sin(\phi_i - \phi_j)\right] - C_s \sin(2\phi_i)$$ |
| 4 & 5 | **Objective function:** $$H = -\sum_{m=1}^{M} \left(\sum_{\substack{i,j \\ i<j}}^{N} \left(-c_{mi}c_{mj}s_i s_j\right) + \sum_{\substack{i,j,k,l \\ i<j<k<l}}^{N} \left(-c_{mi}c_{mj}c_{mk}c_{ml}s_i s_j s_k s_l\right)\right)$$ **Energy function:** |

$$E = C \sum_{m=1}^{M} \left[ \sum_{i,j,i<j}^{N} c_{mi}c_{mj} \cos(\phi_i - \phi_j) \right.$$

$$\left. + \sum_{\substack{i,j,k,l \\ i<j<k<l}}^{N} c_{mi}c_{mj}c_{mk}c_{ml} \cos(\phi_i - \phi_j + \phi_k - \phi_l) + 1 \right]$$

$$- \frac{C_s}{2} \sum_{i=1}^{N} \cos(2\phi_i)$$

**Dynamics:**

$$\frac{d\phi_i}{dt} = C \sum_{m=1}^{M} \left[ \sum_{j=1}^{N} c_{mi}c_{mj} \sin(\phi_i - \phi_j) \right.$$

$$\left. + \sum_{\substack{i \neq j \neq k \neq l \\ j<k<l}}^{N} c_{mi}c_{mj}c_{mk}c_{ml} \sin(\phi_i - \phi_j + \phi_k - \phi_l) \right] - C_s \sin(2\phi_i)$$

*Table 4. Objective functions, corresponding energy expressions, and system dynamics for NAE-K-SAT problems for K=2,3,4, and 5. We note that while the form of the expressions for K=2 and K=3, as well as K=4 and K=5 are similar, the coefficients ($c_{mi}$) are different.*

We now show that the 'energy' functions described in Table 4 decrease with time i.e., $\frac{dE}{dt} \leq 0$ using the example of the energy function for the NAE-4-SAT problem.

$$\frac{\partial E}{\partial \phi_i} = -C \sum_{m=1}^{M} \left[ \sum_{j=1}^{N} c_{mi}c_{mj} \sin(\phi_i - \phi_j) \right.$$

$$\left. + \sum_{\substack{i \neq j \neq k \neq l \\ j<k<l}}^{N} c_{mi}c_{mj}c_{mk}c_{ml} \sin(\phi_i - \phi_j + \phi_k - \phi_l) \right] + C_s \sin(2\phi_i) \quad (3)$$

$$= -\frac{d\phi_i}{dt}$$

Hence,

$$\frac{dE}{dt} = \sum_{i=1}^{N} \frac{\partial E}{\partial \phi_i} \frac{d\phi_i}{dt} = \sum_{i=1}^{N} \left(-\frac{d\phi_i}{dt}\right) \frac{d\phi_i}{dt} = -\sum_{i=1}^{N} \left(\frac{d\phi_i}{dt}\right)^2 \quad (4)$$

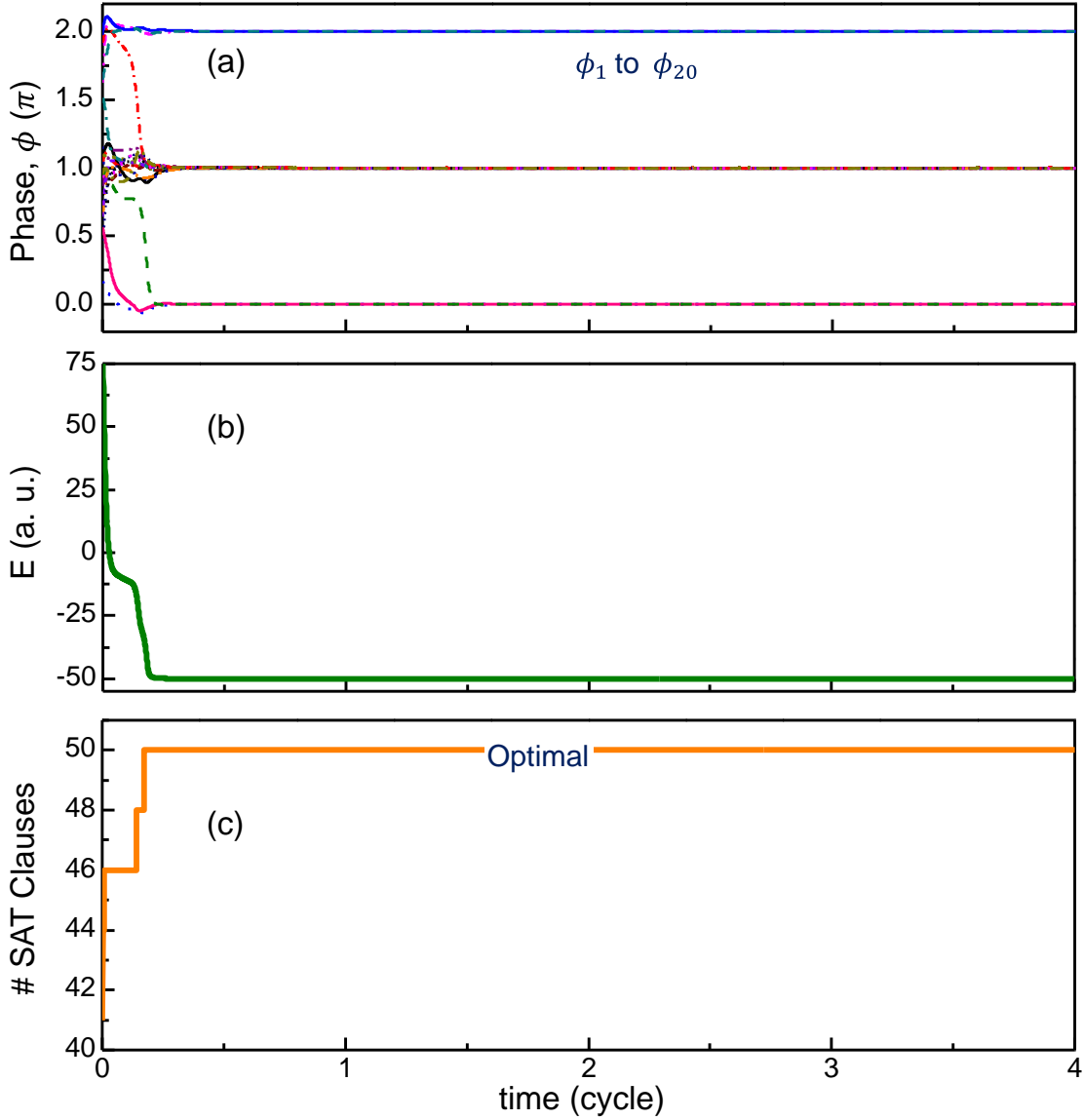

**Figure 1.** Evolution of (a) phases ($\phi$); (b) energy; (c) number of satisfied NAE-4-SAT clauses for an illustrative NAE-4-SAT problem (20 variables and 50 clauses) computed using the proposed dynamical system.

Equation (4) shows that $E$ is a Lyapunov function for the system dynamics formulated for the NAE-4-SAT problem (and the NAE-K-SAT problem in general). Fig. 1 shows an illustrative example of the NAE-4-SAT problem computed using the proposed dynamical system. Details of the simulation used to simulate the illustrative NAE-4-SAT problem are described in Appendix II.

**Max-K-Cut on a Hypergraph:** In the prior section, we exploited the binary nature of the Ising spins (along with higher order interactions among them). We now 'extend' the definition of the 'spin' in order to facilitate the design of computational models for an even broader spectrum of COPs that would benefit from the use of >2 states for each node/spin. To facilitate this, we express the possible states of a spin as $re^{i\theta_k}$, where $r = 1$, and $\theta_k = \frac{2\pi k}{K}$; $k = 1, 2, ... K - 1$. When $K = 2$, the possible states are within {1, -1}, which represents the traditional definition of an Ising spin. In contrast, when $K > 2$, the 'spin' assumes $K$ configurations, represented as complex quantities (e.g., for $K = 3$, the possible states are 1, $e^{i\frac{2\pi(1)}{3}} e^{i\frac{2\pi(2)}{3}}$). While we have utilized this concept for solving combinatorial problems on *graphs* (i.e., problems with quadratic objective functions) [17], here we explore this concept for hypergraphs (that entail higher order interactions) by considering the example of solving the Max-K-Cut of a hypergraph.

Computing the Max-K-Cut on a hypergraph is defined as the challenge of partitioning the nodes of a hypergraph into $K$ partitions in a manner that maximizes the number of hyperedges having nodes that lie in at least two sets created by the partitions [26]. To develop the objective function for the problem, each hyperedge of the graph can be expressed as $h_m = \prod_{i=1}^{N-1} \prod_{j=i+1}^{N} \left(1 - c_{mi}c_{mj}\left(\frac{1-Re\left(s_i s_j^* e^{if(\Delta\theta_{ij})}\right)}{2}\right)\right)$, where $s_j = 1e^{i\theta_j}$; $\theta_j$ can assume any of the following values from $\frac{2\pi k}{K}$; $k = 1, 2, ... K - 1$ enforced by the higher order harmonic injection. $c_{mj} = 1(0)$ if the $j^{th}$ node belongs (does not belong) to the $m^{th}$ hyperedge. We note that the $'i'$ represents the imaginary number $\sqrt{-1}$ whereas $'i'$ refers to the index.

$$f(\Delta\theta_{ij}) = \lim_{\sigma \to 0} \sum_{k=1}^{K-1} \left[ \left((2k-1)\pi - \frac{2k\pi}{K}\right) \cdot e^{-\left(\frac{\left(\Delta\theta_{ij} - \frac{2k\pi}{K}\right)^2}{2\sigma^2}\right)} \right.$$

$$\left. + \left(\frac{2k\pi}{K} - (2k-1)\pi\right) \cdot e^{-\left(\frac{\left(\Delta\theta_{ij} + \frac{2k\pi}{K}\right)^2}{2\sigma^2}\right)} \right] \quad (5)$$

$f(\Delta\theta_{ij})$ is designed such that $Re(s_i s_j^* e^{if(\Delta\theta_{ij})}) = -1(1)$, if the nodes $i$ and $j$ are placed in different (same) sets, and essentially rewards (penalizes) the system in terms of energy, respectively. Additional details about the design and properties of $f(\Delta\theta_{ij})$ have been presented in our prior work [17]. Consequently, if the hyperedge satisfies the criterion for the Max-K-Cut i.e., that the nodes that are connected by it belong to at least two sets, the corresponding assumes $h_m$ assumes a value of 0, else $h_m = 1$. Subsequently, the objective function for the problem, which entails maximizing the number of such hyperedges, can be expressed as minimizing $H$, where,

$$H = \sum_{m=1}^{M} h_m \equiv \sum_{m=1}^{M} \prod_{i=1}^{N-1} \prod_{j=i+1}^{N} \left(1 - c_{mi} c_{mj} \left(\frac{1 - Re(s_i s_j^* e^{if(\Delta\theta_{ij})})}{2}\right)\right) \quad (6)$$

As an example, considering a hypergraph where the maximum number of nodes connected by a hyperedge is 3, the objective function for the Max-K-Cut problem can be expressed as:

$$H = \sum_{\substack{m=1, i \neq j \neq k \\ c_{mi} c_{mj} c_{mk} \neq 0}}^{M} \left(1 - c_{mi} c_{mj} \left(\frac{1 - Re\left(s_i s_j^* e^{if(\Delta \theta_{ij})}\right)}{2}\right)\right) \left(1 - c_{mi} c_{mk} \left(\frac{1 - Re\left(s_i s_k^* e^{if(\Delta \theta_{ik})}\right)}{2}\right)\right) \left(1 - c_{mj} c_{mk} \left(\frac{1 - Re\left(s_j s_k^* e^{if(\Delta \theta_{jk})}\right)}{2}\right)\right) \quad (7)$$

where,

$$f(\Delta \theta_{ij}) = \lim_{\sigma \to 0} \sum_{k=1}^{2} \left[ \left((2k-1)\pi - \frac{2k\pi}{3}\right) \cdot e^{-\left(\frac{\left(\Delta \theta_{ij} - \frac{2k\pi}{3}\right)^2}{2\sigma^2}\right)} + \left(\frac{2k\pi}{3} - (2k-1)\pi\right) \cdot e^{-\left(\frac{\left(\Delta \theta_{ij} + \frac{2k\pi}{3}\right)^2}{2\sigma^2}\right)} \right] \quad (8)$$

For a hypergraph with hyperedges having more than 3 nodes, the objective function entails the use of higher order interactions among the spins.

To formulate a dynamical system for minimizing the above objective function, we express $Re\left(s_i s_j^* e^{if(\Delta \theta_{ij})}\right)$ as $\cos\left(\Delta \theta_{ij} + f(\Delta \theta_{ij})\right)$. Furthermore, we restrict the configuration space of $\theta$ to $\frac{2\pi k}{K}$ where $k = 1, 2, \ldots K - 1$, by injecting the K[th] harmonic (of sufficient strength) which lowers the energy at specific phase points, as described in prior work [17]. The resulting energy function can be described as,

$$E = A \sum_{m=1}^{M} \prod_{i=1}^{N-1} \prod_{j=i+1}^{N} \left(1 - c_{mi}c_{mj}\left(\frac{1 - \cos(\Delta\phi_{ij} + f(\Delta\phi_{ij}))}{2}\right)\right) - \frac{A_s}{K} \sum_{i=1}^{N} \cos(K\phi_i) \quad (9)$$

We note that $\phi$ has been used to express the energy function for the dynamical system instead of $\theta$ which represents the configuration space of the 'extended spin'. The corresponding dynamics for which the function in equation (9) is a Lyapunov function are given by:

$$\frac{d\phi_i}{dt} = -\frac{\partial E}{\partial \phi_i} \quad (10a)$$

$$\frac{d\phi_i}{dt} = \frac{A}{2} \sum_{m=1}^{M} \sum_{j=1, j\neq i}^{N} \left[ c_{mi}c_{mj} \sin(\Delta\phi_{ij} + f(\Delta\phi_{ij})) \frac{h_m}{\left(1 - c_{mi}c_{mj}\left(\frac{1 - \cos(\Delta\phi_{ij} + f(\Delta\phi_{ij}))}{2}\right)\right)} \right] - A_s \sin(K\phi_i) \quad (10b)$$

In the derivation of equation (10b), we exploit the fact that $\frac{\partial f(\Delta\phi_{ij})}{\partial \phi_i} = 0$ [17]. Furthermore, using equation (10a), it can be shown that $\frac{dE}{dt} = -\sum_{i=1}^{N} \left(\frac{d\phi_i}{dt}\right)^2 \leq 0$ (similar to equation (4)).

We now evaluate our proposed model on a representative hypergraph. We consider a hypergraph where each hyperedge has 3 vertices. The corresponding dynamics for this case can then be written as,

$$\begin{aligned}
\frac{d\phi_i}{dt} = \frac{A}{2} \sum_{m=1}^{M} &\left[ c_{mi}c_{mj} \sin\left(\Delta\phi_{ij} + f(\Delta\phi_{ij})\right) \left(1 - c_{mi}c_{mk}\left(\frac{1 - \cos(\Delta\phi_{ik} + f(\Delta\phi_{ik}))}{2}\right)\right) \left(1 \right.\right. \\
&\left. - c_{mj}c_{mk}\left(\frac{1 - \cos\left(\Delta\phi_{jk} + f(\Delta\phi_{jk})\right)}{2}\right)\right) + c_{mi}c_{mk}\sin\left(\Delta\phi_{ik} \right. \\
&\left. + f(\Delta\phi_{ik})\right)\left(1 - c_{mi}c_{mj}\left(\frac{1 - \cos(\Delta\phi_{ij} + f(\Delta\phi_{ij}))}{2}\right)\right)\left(1 \right. \\
&\left.\left. - c_{mj}c_{mk}\left(\frac{1 - \cos(\Delta\phi_{jk} + f(\Delta\phi_{jk}))}{2}\right)\right)\right] - A_s \sin(K\phi_i)
\end{aligned} \quad (11)$$

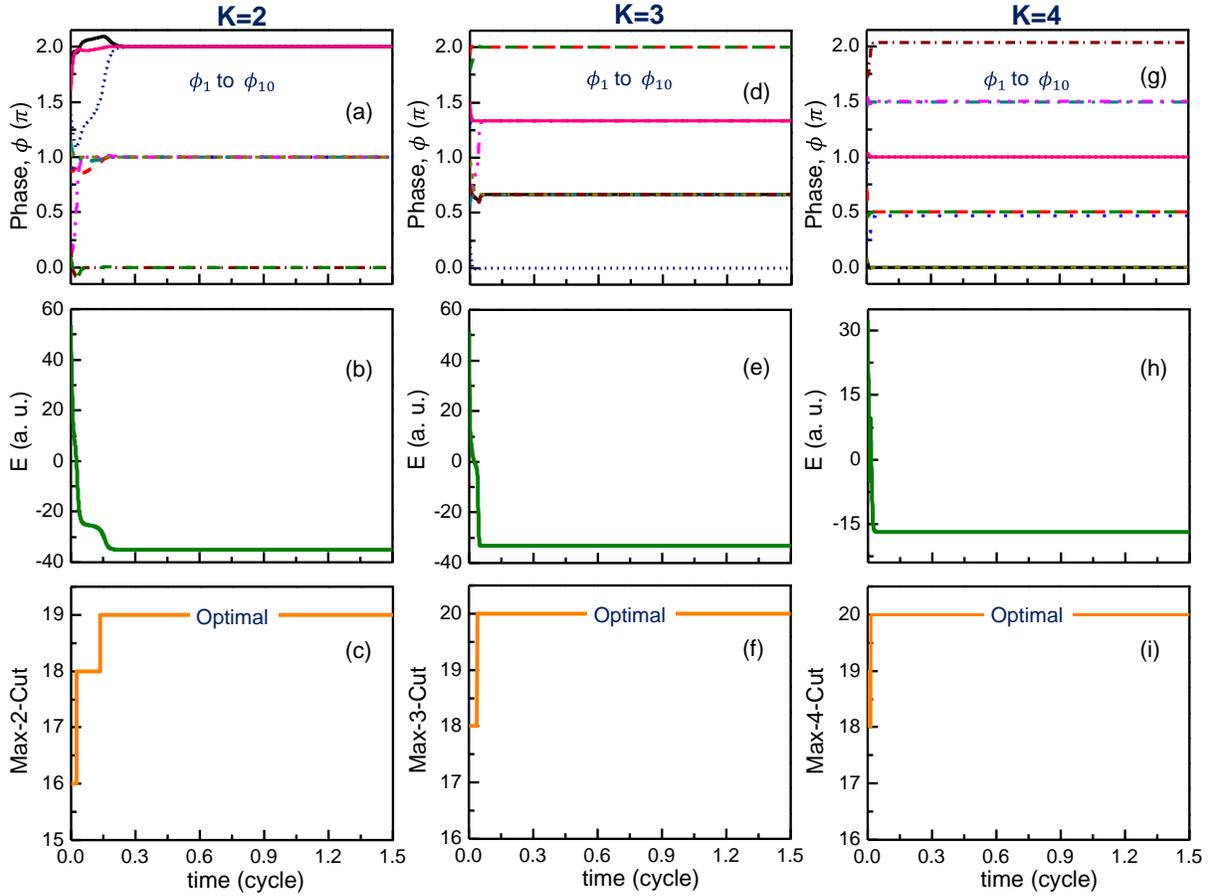

**Figure 2.** Max-K-Cut (K=2,3,4) solutions computed using the proposed dynamical system for an illustrative hypergraph. Evolution of phases ($\phi$), energy and the Max-K-Cut solution, respectively for (a-c) K=2; (d-f) K=3; (g-i) K=4.

Fig. 2 shows the computed Max-K-Cut (for K=2, 3, and 4) for a hypergraph instance (with 10 nodes, and 20 hyperedges). The illustrative problem has a maximum of 4 nodes per hyperedge. Details of the simulation used to simulate the illustrative Max-K-Cut problem are described in Appendix II.

| No. Partitions / No. of spin configurations \ Interaction order | Quadratic (=2) Representation: Graph | Higher Order (>2) Representation: Hypergraph |
|---|---|---|
| **Binary** | -MaxCut<br>-Spin Glass<br>Bashar *et al.*, arXiv (2022) [27]<br>Bashar et al., JxCDC (2020) [3] | -Boolean SAT<br>-NAE-SAT /Set-Splitting<br>This Work |
| **Multi-valued** | -Max-K-Cut<br>-Graph coloring<br>-TSP<br>Mallick *et al.*, Phys. Rev App (2022) [17] | -Hypergraph Max-K-Cut<br>-Hypergraph colorig<br>This Work |

**Figure 3.** Proposed classification of COPs based on the number of states/configurations for the nodes, and the nature of interaction among them. The proposed work enables the direct development of computational models for COPs that entail higher order interactions.

**Discussion**

In this work, we develop computational models, inspired by dynamical systems, that model higher order interactions (beyond quadratic/pairwise) among Ising spins. Our approach enables the direct formulation of analog computing models for many COPs that entail such interactions without the need for problem decomposition. Furthermore, using the combination of higher order interactions along with 'expanding' the number of 'spin' states to greater than 2, we can directly map and solve an even broader class of problems on hypergraphs. This has been summarized in Fig. 3. Consequently, in the context of the broader effort focused on developing dynamical system-inspired models for solving hard COPs, this work expands on the potential of physics-inspired solvers to accelerate COPs.

**Appendix I**

Here, we develop the formulation of the objective function and the corresponding dynamical system for the NAE-5-SAT (K is an odd number) problem. An NAE-5-SAT clause can be represented as,

$$C = (x_i \oplus x_j) \vee (x_i \oplus x_k) \vee (x_i \oplus x_l) \vee (x_i \oplus x_m) \vee (x_j \oplus x_k) \vee (x_j \oplus x_l)$$
$$\vee (x_j \oplus x_m) \vee (x_k \oplus x_l) \vee (x_k \oplus x_m) \vee (x_l \oplus x_m) \quad \text{(A1.1)}$$

In terms of Ising spins, the complement of $C$ can be written as,

$$\left(\frac{1+s_i s_j}{2}\right)\left(\frac{1+s_i s_k}{2}\right)\left(\frac{1+s_i s_l}{2}\right)\left(\frac{1+s_i s_m}{2}\right)\left(\frac{1+s_j s_k}{2}\right)\left(\frac{1+s_j s_l}{2}\right)\left(\frac{1+s_j s_m}{2}\right)$$

$$\left(\frac{1+s_k s_l}{2}\right)\left(\frac{1+s_k s_m}{2}\right)\left(\frac{1+s_l s_m}{2}\right)$$

$$= \frac{1}{2^4}\big(1 + s_i s_j + s_i s_k + s_i s_l + s_i s_m + s_j s_k + s_j s_l + s_j s_m + s_k s_l \quad \text{(A1.2)}$$

$$+ s_k s_m + s_l s_m + s_i s_j s_k s_l + s_i s_j s_k s_m + s_i s_j s_l s_m + s_i s_k s_l s_m$$

$$+ s_j s_k s_l s_m\big)$$

Thus, the objective function for an NAE-5-SAT problem with M clauses can be written as,

$$H = -\sum_{m=1}^{M}\left(\sum_{\substack{i,j \\ i<j}}^{N}(-c_{mi}c_{mj}s_i s_j) + \sum_{\substack{i,j,k,l \\ i<j<k<l}}^{N}(-c_{mi}c_{mj}c_{mk}c_{ml}s_i s_j s_k s_l)\right) \quad \text{(A1.3)}$$

Where, $c_{mi} = -1 \,(+1)$ if the $i^{th}$ variable appears in inverted (normal) form in the $m^{th}$ clause; $c_{mi} = 0$ if the $i^{th}$ variable is absent in the $m^{th}$ clause.

Using the approach described in the main text, the corresponding Lyapunov function and the system dynamics can be formulated as,

Energy:

$$E = C \sum_{m=1}^{M} \left[ \sum_{i,j,i<j}^{N} c_{mi} c_{mj} \cos(\phi_i - \phi_j) \right.$$

$$\left. + \sum_{\substack{i,j,k,l \\ i<j<k<l}}^{N} c_{mi} c_{mj} c_{mk} c_{ml} \cos(\phi_i - \phi_j + \phi_k - \phi_l) + 1 \right]$$

$$- \frac{C_s}{2} \sum_{i=1}^{N} \cos(2\phi_i)$$

(A1.4)

Dynamics:

$$\frac{d\phi_i}{dt} = C \sum_{m=1}^{M} \left[ \sum_{j=1}^{N} c_{mi} c_{mj} \sin(\phi_i - \phi_j) \right.$$

$$\left. + \sum_{\substack{i \neq j \neq k \neq l \\ j<k<l}}^{N} c_{mi} c_{mj} c_{mk} c_{ml} \sin(\phi_i - \phi_j + \phi_k - \phi_l) \right] - C_s \sin(2\phi_i)$$

(A1.5)

Equations (A1.3), (A1.4), and (A1.5) are also shown in Table 4 in the main text.

**Appendix II**

Here, we describe the simulation approach used to simulate the NAE-4-SAT problem (Fig. 1, main text) and the hypergraph Max-K-Cut problem (Fig. 2, main text). We solve the dynamics using a stochastic differential equation (SDE) solver implemented in MATLAB; details of its implementation have been described in our previous work [28]. The SDE solver incorporates noise that helps escape local minima in the phase space.

Values of $C$ and $C_s$ used in the simulation of the NAE-4-SAT are:

| Problem Solved | $C$ | $C_s$ |
|---|---|---|
| NAE-4-SAT | $\frac{10}{8}$ | 5 |

Values of $A$ and $A_s$ used in the simulation of the Max-K-Cut are:

| Problem Solved | $A$ | $A_s$ |
|---|---|---|
| Hypergraph Max-2-Cut | 15 | 10 |
| Hypergraph Max-3-Cut | 15 | 10 |
| Hypergraph Max-4-Cut | 10 | 10 |

**Acknowledgment**

This work was supported by NSF ASCENT grant (No. 2132918). We would like to thank Professor Zongli Lin from the University of Virginia for insightful discussions.

6. Marandi, A., Wang, Z., Takata, K., Byer, R. L. & Yamamoto, Y. Network of time-multiplexed optical parametric oscillators as a coherent Ising machine. Nat. Photonics **8**, 937–942 (2014).

7. McMahon, P.L., Marandi, A., Haribara, Y., Hamerly, R., Langrock, C., Tamate, S., Inagaki, T., Takesue, H., Utsunomiya, S., Aihara, K. & Byer, R.L. A fully programmable 100-spin coherent Ising machine with all-to-all connections. Science **354**, 614-617 (2016).

8. Böhm, F., Verschaffelt, G. and Van der Sande, G. A poor man's coherent Ising machine based on opto-electronic feedback systems for solving optimization problems. Nature communications **10**, 1-9 (2019).

9. Mohseni, N., McMahon, P.L. and Byrnes, T. Ising machines as hardware solvers of combinatorial optimization problems. Nature Reviews Physics **4**, 363-379 (2022).

10. Vadlamani, S.K., Xiao, T.P., and Yablonovitch, E. Physics successfully implements Lagrange multiplier optimization, Proceedings of the National Academy of Sciences **117**, 26639-26650 (2020).

11. Lucas, A. Ising formulations of many NP problems. Frontiers in physics **5** (2014).

12. Terada, K., Oku, D., Kanamaru, S., Tanaka, S., Hayashi, M., Yamaoka, M., Yanagisawa, M. and Togawa, N. An Ising model mapping to solve rectangle packing problem. In 2018 International Symposium on VLSI Design, Automation and Test (VLSI-DAT), 1-4, IEEE (2018).
23

21. Mandal, A., Roy, A., Upadhyay, S. and Ushijima-Mwesigwa, H. Compressed quadratization of higher order binary optimization problems. In Proceedings of the 17th ACM International Conference on Computing Frontiers, 126-131 (2020).

22. Choi, V. Adiabatic quantum algorithms for the NP-complete Maximum-Weight Independent set, Exact Cover and 3SAT problems. arXiv preprint arXiv:1004.2226 (2010).

23. Sanyal, S. Reduction from SAT to 3SAT. Available at: https://cse.iitkgp.ac.in/~palash/2018AlgoDesignAnalysis/SAT-3SAT.pdf

24. Abel, S.A. and Nutricati, L.A. Ising Machines for Diophantine Problems in Physics. Fortschritte der Physik, 2200114 (2022).

25. Wang T., and Roychowdhury, J. Oscillator-based Ising machine. arXiv preprint (2017). At: https://arxiv.org/pdf/1903.07163.pdf

26. Ageev, A. A. and Sviridenko, M. I. An approximation algorithm for hypergraph max k-cut with given sizes of parts. In European Symposium on Algorithms, 32-41, Springer, Berlin, Heidelberg (2000).

27. Bashar, M.K., Mallick, A., Ghosh, A.W. and Shukla, N. Dynamical system-based computational models for solving combinatorial optimization on hypergraphs. arXiv preprint arXiv:2207.09618 (2022).

28. Bashar, M.K., Lin, Z. and Shukla, N. Formulating Oscillator-Inspired Dynamical Systems to Solve Boolean Satisfiability. arXiv preprint arXiv:2209.07571 (2022).
25